\documentclass[10pt]{amsart}
\usepackage{geometry}                
\usepackage{graphicx}
\usepackage{amssymb}
\usepackage{epstopdf}
\newtheorem{theorem}{Theorem}[section]
\newtheorem{lemma}[theorem]{Lemma}
\newtheorem{corollary}[theorem]{Corollary}

\theoremstyle{definition}

\def\eproof{$\Box$ \medskip}
\def\sech{\mbox{sech}}
\newcommand{\Real}{\mathbb R}

\newcommand{\Hyp}{\mathbb{H}}
\newcommand{\Sph}{\mathbb{S}}

\newcommand{\Cplex}{\mathbb{C}}
\newcommand{\Vol}{\mbox{Vol}}

\title[Moments of hitting function]{Moments of the boundary hitting function for the geodesic flow on a hyperbolic manifold}
\author{Martin Bridgeman}
\address{Boston College}
\author{Ser Peow Tan}
\address{National University of Singapore}
\date{}                                           

\thanks{Tan was partially supported by the National University
of Singapore academic research grant R-146-000-156-112.}

\begin{document}

\maketitle

\begin{abstract}    
In this paper we consider finite volume hyperbolic manifolds $X$ with non-empty totally geodesic boundary. We consider the distribution of the times for the geodesic flow to hit the boundary and derive a formula for the moments of the associated random variable in terms of the orthospectrum. We show that the the first two moments correspond to two cases of known identities for the orthospectrum. We further obtain an explicit formula in terms of the trilogarithm functions for the average time for the geodesic flow to hit the boundary in the surface case, using the third moment.
\end{abstract}


\section{Introduction}
Let $X$ be a finite volume hyperbolic manifold with non-empty totally geodesic boundary  $\partial X$. An {\em orthogeodesic} for $M$ is a geodesic arc with endpoints perpendicular to $\partial X$. These were first introduced by Basmajian in \cite{Bas93} in the study of totally geodesic submanifolds. We denote by $O_X = \{\alpha_{i}\}$  the collection of orthogeodesics of $X$ and let $l_i$ be the length of $\alpha_{i}$. We note that $O_X$ is countable as the elements correspond to a subset of the collection of closed geodesics of the  double of $X$ along its boundary.  We call the set $L_X = \{l_{i}\}$ (with multiplicities) the {\em orthospectrum.}

In \cite{Bas93}, Basmajian derived the following boundary orthospectrum identity;
\begin{equation} \mbox{Vol}(\partial X) = 2 \sum_{l \in L_{X}} V_{n-1}\left(\log\left(\coth{\frac{l}{2}}\right)\right)
\label{basmajian}
\end{equation}
where $V_n(r)$ is the volume of the ball of radius $r$ in $\Hyp^n$.
The identity comes from considering the universal cover $\tilde{X} \subseteq \Hyp^n$ of $X$. Then $\partial\tilde{X}$ is a countable collection of disjoint hyperbolic hyperplanes  which are the lifts of the boundary components of $\partial X$. For each component $C$ of $\partial\tilde{X}$, we orthogonally project each of the other components of $\partial\tilde{X}$ onto $C$ to obtain a collection of disjoint disks on each component $C$ of $\partial\tilde{X}$. These disks form an equivariant family of disks that are full measure in $\partial\tilde{X}$. They descend to a family of disjoint disks in $\partial X$ of full measure. As each orthogeodesic lifts to a perpendicular between two components of $\partial\tilde{X}$, each orthogeodesic corresponds to two disks (one at each end) in the family of disks in $\partial X$  and this  gives the above identity.

Using a decomposition of the unit tangent bundle, Bridgeman-Kahn (see \cite{BK10})  derive the identity;
\begin{equation}
 \Vol(T_1(X)) = \sum_{l \in L_{X}} H_n(l)
 \label{bridgeman}
 \end{equation}
 where $H_n$ is some smooth function depending only on the dimension $n$. As $\Vol(T_1(X)) = \Vol(X).V_{n-1}$ where $V_{n-1}$ is the volume of the unit sphere in $\Real^n$, the above identity can also be written as
$$ \Vol(X) = \sum_{l \in L_{X}} \overline{H}_n(l) $$
where $\ \overline {H}_n(l) = H_n(l)/V_{n-1}$.\newline  In the specific case of dimension two, the function $ H_n$ is given in terms of the Rogers dilogarithm (see \cite{B11}). In the papers \cite{Cal10a,Cal10b} Calegari gives an alternative derivation of the identity  in equation \ref{bridgeman}.

The motivation for this paper was  to connect the above two identities in a natural framework. The connection is that they are the first two moments of the Liouville measure. A second motivation was to compute the average time it takes to hit the boundary under the geodesic flow. This can be put into the same framework and it turns out that consideration of the third moment gives a formula for the average time it takes to hit the boundary of $M$ under the geodesic flow in terms of the orthospectrum. It is conceivable that higher moments encode other important geometric invariants of the manifold.

\section{Moments of Liouville measure}
We let $G(\Hyp^n)$ be the space of oriented geodesics in $\Hyp^n$.  By identifying a geodesic with its endpoints on the sphere at infinity,  the space $G(\Hyp^n) \simeq (\Sph_\infty^{n-1}\times \Sph_\infty^{n-1} - \mbox{Diagonal})$. The {\em Liouville measure} $\mu$ on $G(\Hyp^n)$ is a Mobius invariant measure. In the upper half space model, we identify a geodesic with its endpoints $(x,y) \in  \overline{\Real^{n-1}}\times\overline{\Real^{n-1}}$. Then the Liouville measure $\mu$ has the form  $$d\mu_{(x,y)} = \frac{2dV_xdV_y}{|x-y|^{2n-2}}$$
where $dV_x = dx_1dx_2\ldots dx_{n-1}$, for $x = (x_1, x_2,\ldots,x_{n-1}) \in \Real^{n-1}$.

If $X$ is a hyperbolic $n$-manifold with totally geodesic boundary, we identify $\tilde{X}$, the universal cover of $X$ as a subset of $\Hyp^n$ and  $\Gamma  \subseteq \mbox{Isom}^{+}(\Hyp^n)$ such that $X = \tilde{X}/\Gamma$. Then $G(\tilde{X}) \subseteq G(\Hyp^n)$  is the set of geodesics intersecting $\tilde{X}$. We define $G(X) = G(\tilde{X})/\Gamma$, the space of geodesics in $X$. Then by invariance of  Liouville measure $\mu$ descends to a measure on $G(X)$ which we also call $\mu$.

We define the  measurable function $L:G(X) \rightarrow [0,\infty]$  by  $L(g) := \mbox{Length}(g)$, where the length is taken in $X$. This is the  hitting length function for $X$.  As the limit set $L_\Gamma$ has measure zero, almost every geodesic hits the boundary  of $X$ and therefore for almost every geodesic $g$, $L(g)$ is finite and $g$ is a proper geodesic arc.

We define the pushforward measure $M := L_*(\mu)$ on the real line. This measure is the  distribution of lengths of geodesics in $X$. We define its k-th moment to be
$$M_k(X) = M(x^k) = \int_0^\infty x^k dM = \int_{G(X)} L^k(g) \ d\mu.$$
In general, the moments of a random variable give a set of measurements that describe distributional properties of the random variable such as the average value and variance. Using a decomposition of $G(X)$, we show that the moments $M_k(X)$ have formulae that extend the identities in equations \ref{basmajian},\ref{bridgeman}.

\medskip

The main result of the paper is the following:

\medskip
\begin{theorem}(Main Theorem)
{\em There exists smooth functions $F_{n,k}: \Real_+ \rightarrow \Real_+$ and constants $C_n > 0$  such that if $X$ is a compact hyperbolic n-manifold with totally geodesic boundary $\partial X \neq \emptyset$, then
\begin{enumerate}
\item The moment $M_k(X)$ satisfies
$$M_k(X) = \sum_{l \in L_X} F_{n,k}(l)$$
\item $M_{0}(X) = C_n.\Vol(\partial X)$ and  the identity for $M_{0}(X)$ is the  identity in equation \ref{basmajian}.
\item $M_{1}(X) = \Vol(T_1(X))$ and the identity for $M_1(X)$ is the identity in equation \ref{bridgeman}.
\item  $M_2(X) = 2\Vol(T_1(X))A(X)$ where $A(X)$ is the average time for a vector in $T_1(X)$ to hit the boundary  under geodesic flow. Therefore by the identity for $M_2(X)$
$$A(X) = \frac{1}{2\Vol(T_1(X))}\sum_{l \in L_X} F_{n,2}(l) = \sum_{l \in L_X} G_{n}(l).$$
\end{enumerate}
}
\end{theorem}

In the surface case we obtain an explicit formula for the function $G_2$ and hence $A(X)$ in terms of polylogarithms. Furthermore, besides compact surfaces obtained as quotients of Fuchsian groups, the identity holds more generally for  finite area surfaces, which we describe next.

If $S$ is a finite area surface with totally geodesic boundary $\partial S \neq \emptyset$, then the boundary components are either closed geodesics or bi-infinite geodesics with cuspidal endpoints. We define a {\em boundary cusp} of $S$ to be an ideal vertex of $\partial S$. We let $C_S$ be the number of boundary cusps of $S$. Then we have the following explicit formula for $A(S)$:

\begin{theorem}
Let $S$ be a finite area hyperbolic surface with non-empty totally geodesic boundary. Then
$$A(S) = \frac{1}{8\pi^2|\chi(S)|}\left(\sum_{l \in L_S} F\left(\frac{1}{\cosh^2(l/2)}\right) + 6\zeta(3)C_S\right)$$
where
\begin{eqnarray*}F(a) &=&  -12\zeta(3)-\frac{4\pi^2}{3}\log(1-a)+6\log^2(1-a)\log(a)-4\log(1-a)\log^2(a)\\
&& \qquad -8\log\left(\frac{a^2}{1-a}\right)Li_2(a)
+24Li_3(a)+12Li_3(1-a),
\end{eqnarray*}
$Li_k(x)$ is the $k^{th}-$polylogarithm function, and $\zeta$ is the Riemann $\zeta-$function.
\label{avehit}
\end{theorem}

\section{A natural fibering}
We have the natural fiber bundle $p:T_1(\Hyp^n) \rightarrow G(\Hyp^n)$ such that $v$ is tangent to the oriented geodesic $p(v)$. Let $\Omega$ be the volume measure on $T_1(\Hyp^n)$.  We parametrize $T_1(\Hyp^n)$ as follows: We first choose a basepoint $b_g$ on each geodesic $g$ (say by taking the point closest to a fixed point p). Then a vector $v \in T_1(\Hyp^n)$ is given by a triple $(x,y,l) \in \overline{\Real^{n-1}}\times\overline{\Real^{n-1}}\times\Real$ where $v$ is tangent to the geodesic $g$ with endpoints $x,y$ (from $x$ to $y$), and $l$ is the signed length along $g$ from the basepoint $b_g$. In terms of this parametrization the volume form $\Omega$ on $T_1(\Hyp^n)$ is
$$d\Omega_{v} = \frac{2dV_xdV_ydl}{|x-y|^{2n-2}} = d\mu_{(x,y)}dl_{v}$$
where $\mu$ is the Liouville measure (see \cite{Nic89}).
If $X$ is a hyperbolic $n$-manifold with totally geodesic boundary and $\tilde{X}$ the universal cover of $X$, then the fiber bundle $p$ restricts to $T_1(\tilde{X})$ to give the equivariant map
$p:T_1(\tilde{X}) \rightarrow G(\tilde{X})$ which   descends to a map $\overline{p}:T_1(X) \rightarrow G(X)$.
We have the function $\overline{L}: T_1(X) \rightarrow [0,\infty]$ given by $\overline{L} = L\circ \overline{p}$.  The measure  $N = \overline{L}_*\Omega$ was introduced in \cite{B11} in order to derive the surface case of the identity \ref{bridgeman}. The measures $M,N$ have a simple relation which we now describe.

\medskip
For $\phi:[0,\infty)\rightarrow\Real$ a smooth function with compact support then
$$N(\phi) = \Omega(\phi\circ \overline{L}) = \int_{v \in T_1(X)} \phi(\overline{L}(v))d\Omega_v = \int_{g \in G(X)} \int_{v \in p^{-1}(g)}\phi(L(g))d\mu_gdl_v$$
$$ = \int_{g \in G(X)}\phi(L(g))\left( \int_{v \in p^{-1}(g)}dl_v\right)d\mu_g =   \int_{g \in G(X)}\phi(L(g))L(g) d\mu_g = M(x\phi)$$
It follows that the measures $M, N$  satisfy $ dN = xdM.$ We define the moments of $N$ to be $N_k(X) = N(x^k)$. Then
\begin{equation}
N_k(X) = M_{k+1}(X).
\label{2moments}
\end{equation}
Also as $N_0(X) = \Vol(T_1(X))$ it follows that $M_1(X) = \Vol(T_1(X))$.

Because of the above relation between the measures $M,N$, the results described in this paper can be given in terms of either.  For the most part we will give our results in terms of  the measure $M$ and its moments, but when it is more natural to do so, we will consider the measure $N$.

We let $A(X)$ be the average time for a vector in $T_1(X)$ to hit the boundary under geodesic flow. Then as the time for $v$ and $-v$ sum to $\overline{L}(v)$,  $A(X)$ is half the average of the function  $\overline{L}$. Then $A(X)$ is given by the first moment of measure $N$ and as $N_1(X) = M_2(X)$,  we
 have the formula
$$A(X) = \frac{1}{2}\frac{M_2(X)}{\Vol(T_1(X))}$$
proving part 4) of the Main Theorem.

\section{Moments are finite}
Before we derive summation formulae for the moments $M_k(X)$, we need to first show that they are finite. The proof that $M_0(X), M_1(X)$ are finite will follow from explicit calculation, in particular, from the last section $M_1(X) = \Vol(T_1(X))$ which is finite by assumption. By equation \ref{2moments}, $M_{k}(X) = N_{k-1}(X)$, and therefore we need only show that $N_k(X)$ is finite for $k \geq 1$. To prove this,  we show that the measure $N$ on the real line decays exponentially, i.e. there exist positive constants $a, C$ such that $dN \leq Ce^{-at}dt$ for $t$ large. Then $N_k(X)$ is finite as the measure $x^ke^{-at}dt$ is finite.

We first recall some background on Kleinian groups (see \cite{Mas87} for details). A {\em Kleinian group} $\Gamma$ is a discrete subgroup of the isometries of $\Hyp^n$. The {\em limit set} $L_\Gamma = \overline{\Gamma x}\cap \Sph^{n-1}_\infty$ is the accumulation set of an orbit of a point $x$ on the boundary. It is  easy to show that $L_\Gamma$ is independent of $x$. The {\em convex hull} of $\Gamma$, denoted $H(\Gamma)$, is the smallest convex set containing all geodesics with endpoints in $L_\Gamma$. As $H(\Gamma)$ is invariant under $\Gamma$, the {\em convex core} is defined to be $C(\Gamma) = H(\Gamma)/\Gamma$. A Kleinian group is {\em convex cocompact} if $C(\Gamma)$ is compact. Also a group is {\em geometrically finite} if $N_\epsilon(C(\Gamma))$, the $\epsilon$ -neighborhood of the core, is finite volume.

Let $\Gamma$ be a convex cocompact  Kleinian group with $N = \Hyp^n/\Gamma$ and $X = H(\Gamma)/\Gamma$ its convex core. We let $\delta(\Gamma)$ be the Hausdorff dimension of the limit set $L_\Gamma$. Let $g_t$ be the geodesic flow on $T_1(N)$ and define
$$B(t) = \{ v\in T_1(X)\ | \ g_t(v) \in T_1(X)\} = g_t(T_1(X))\cap T_1(X) $$
The set $B(t)$ is the set of tangent vectors that remain in the convex core under time t flow.
We now use  a standard counting argument on orbits to bound the volume of the set $B(t)$ (see \cite{Nic89} for background).

\begin{lemma}
Given  $\Gamma$ a convex cocompact Kleinian group, then there exists constants $A, T$ such that
$$\Vol(B(t)) \leq Ae^{-(n-1 -\delta(\Gamma))t}$$
for $t > T$.
In particular if $L_\Gamma \neq \Sph^{n-1}_\infty$ then $\Vol(B(t))$ is exponentially decaying.
\end{lemma}
{\bf Proof:}
We take $0 \in H(\Gamma)$ and consider its orbits under $\Gamma$. Then we let
$$O(r) = \{\gamma \in \Gamma\ | d(0,\gamma(0)) < r\}\qquad \mbox{and} \qquad N(r) = \#O(r).$$
By Sullivan (see \cite{Sul79}), there exists  constants $A, r_0$ such that $N(r) \leq Ae^{\delta(\Gamma)r}$ for $r > r_0$. We let $D$ be the diameter of $X$. Given a unit  tangent vector $v$, we denote its basepoint by $b(v)$. We let $\tilde{B}(t)$ be a lift of $B(t)$ with basepoints within a distance  $D$ of $0$.  For $v \in \tilde{B}(t)$ then by the trangle inequality,  $g_t(v)$ has basepoint $b(g_t(v))$ such that $t-D < d(b(g_t(v)), 0) <t+D$. Therefore $b(g_t(v))$ has nearest orbit $\gamma(0)$ such that $b(g_t(v)) \in B(\gamma(0),D)$, the ball of radius $D$ around $\gamma(0)$. Also we have that  $t-2D < d(\gamma(0),0) <t+2D$. Therefore
$$U(t) = \bigcup_{O(t-2D)^c\cap O(t+2D)} B(\gamma(0), D)$$

Let  $x$ be a distance at most $D$ from $0$. We want to bound the visual measure of $U(t)$ from $x$. Each $B(\gamma(0),D) \in U(t)$ is a distance between $t-3D,t+3D$ from $x$. Let  $T_1 > 3D$, and restrict to   $t > T_1$.  We radially project each $B(\gamma(0),D) \in U(t)$ onto the  $S(x,t)$, the sphere of radius $t$ in hyperbolic space. We label the projection $P(\gamma(0),D)$.  Then  the area of each projection is bounded above by a constant $C_1>0$. Therefore
$$Vis_x(U(t)) \leq \frac{C_1N(t+2D)}{\Vol(S(x,t))}.$$
We have that $\Vol(S(x,t)) = S_n\sinh^{n-1}(t)$ where $S_n$ is the volume of the standard Euclidean sphere of dimension $(n-1)$. Thus for $t > T_1$, $\Vol(S(x,t))\geq  L_ne^{(n-1)t}$ for some $L_n > 0$.

Therefore for
$$Vis_x(U(t)) \leq \frac{C_1.N(t+2D)}{L_ne^{(n-1)t}} \leq \frac{C_1.Ae^{\delta(\Gamma)(t+2D)}}{L_ne^{(n-1)t}} \leq Ce^{-((n-1)-\delta(\Gamma))t}$$
for some constant $C$.

In order to obtain the bound on $\Vol(B(t))$ we integrate the visual measure of $U(t)$ over $\tilde{B}(t)$ gives
$$\Vol(B(t)) = \Vol(\tilde{B}(t)) \leq \Vol(X).Ce^{-((n-1)-\delta(\Gamma))t}$$
for $t > T$, giving our result for $\Gamma$ convex cocompact.
\eproof

\begin{corollary}
If $X$ is a compact hyperbolic manifold with non-empty totally geodesic boundary then the moments $M_k(X)$ are finite for $k \geq 1$.
\end{corollary}
{\bf Proof:}
Let
$$E(t) = \{ v \in T_1(X) \  |\   \overline{L}(v) \in [t,t+1)\}.$$
Then from the above lemma $E(t) \subset B(t/2)\cup B(-t/2)$. Therefore there are constants $a, K > 0$ such that $Vol(E(t)) \leq Ke^{-at}$ for $t > 2T_0$.
Therefore $$N_{k}(X)  \leq \sum_{n=0}^\infty (n+1)^{k}. Vol(E(n))$$
is finite for $k \geq 0$ by comparison with the series $\sum n^ke^{-an}.$ Therefore $M_k(X)$ is finite for $k \geq 1$.
\eproof

 \section{Decomposition of the space of geodesics}
We let $G(\Hyp^n)$ be the space of oriented geodesics in $\Hyp^n$. By identifying a geodesic with its endpoints, the Liouville measure $\mu$ on $G(\Hyp^n)$  is given by
$$d\mu = \frac{2dV_xdV_y}{|x-y|^{2n-2}}$$

If $X$ is a hyperbolic $n$-manifold with totally geodesic boundary then we let $\tilde{X}$ be the universal cover of $X$ as a subset of $\Hyp^n$ and  $\Gamma  \subseteq \mbox{Isom}^{+}(\Hyp^n)$ such that $X = \tilde{X}/\Gamma$. We let $G(\tilde{X}) \subseteq G(\Hyp^n)$ be the set of geodesics intersecting $\tilde{X}$. We define $G(X) = G(\tilde{X})/\Gamma$, the space of geodesics in $X$. Then by invariance of  Liouville measure $\mu$ descends to a measure on $G(X)$ which we also call $\mu$.

The space $G(X)$ has a simple (full measure) decomposition via orthogeodesics. For $\alpha$ an orthogeodesic we define
$$F_\alpha =\{ g \in G(X) \ | \ g \mbox{ is homotopic rel $\partial X$ to } \alpha\}.$$
The set of orthogeodesics is countable, so we index our orthogeodesics $ O_X =\{\alpha_i\}$ and the sets $F_i = F_{\alpha_i}$.
As the limit set $L_\Gamma$ has zero measure, $\bigcup F_i$ gives a (full measure) partition of $G(X)$ with respect to $\mu$. We note that the sets $F_i \subseteq \Sph_\infty^{n-1}\times \Sph_\infty^{n-1}$ are of the form $(D_1\times D_2)\cup(D_2\times D_1)$ where $D_1, D_2$ are disjoint round disks. Thus  $G(X)$  decomposes (up to full measure) into a countable collection of  elementary pieces indexed by the orthogeodesics.

Then
$$M_k(X) = M(x^k) = \int_0^\infty x^k dM = \int_{G(X)} L^k(g)\  d\mu  = \sum_{i} \int_{F_i} L^k(g)\  d\mu$$
where  $L:G(X) \rightarrow \Real$ is the length of a geodesic.

We now lift $F_i$ to a set $\tilde{F}_i \subset G(\Hyp^n)$.  In the upper half-space model, we choose two planes $P, Q$ with orthogonal distance $l_i = Length(\alpha_i)$. Then $F_i$ lifts to the set of geodesics  which intersect $P$ and $Q$. In particular we can  take $P, Q$ to have boundary circles centered at $0$ with radii $1, e^{l_i}$ respectively. Lifting the function $L$ to $\tilde{F}_i$ we see that it only depends on the endpoints $(x, y)$ and the ortholength $l_i$. We denote it by  $L(x,y, l_i)$.
$$ \int_{F_i} L^k(g) d\mu = \int_{\tilde{F}_i} \frac{2L^{k}(x,y) dV_xdV_y}{|x-y|^{2n-2}}$$
Integrating over the length parameter (in both directions) we obtain
$$\int_{F_i} L^k(g) d\mu = \int_{|x| <1} \int_{|y| > e^{l_i}} \frac{4L(x,y, l_i)^{k} dV_xdV_y}{|x-y|^{2n-2}} = F_{n,k}(l_i).$$
Therefore we have that
$$M_{k}(X) = \sum_i F_{n,k}(l_i)$$
 and $F_{n,k}$ is given by the integral formula
 $$F_{n,k}(t) = \int_{|x| <1} \int_{|y| > e^{t}} \frac{4L(x,y,t)^{k} dV_xdV_y}{|x-y|^{2n-2}} $$
 This gives the summation formula for the moment
 $$M_k(X) = \sum_{l \in L_X} F_{n,k}(l)$$
and proves part 1) of the Main Theorem.

By equation \ref{2moments}, $M_1(X) = N_0(X) = \Vol(T_1(X))$. Thus letting $E_i \subseteq T_1(X)$ given by $E_i = \overline{p}^{-1}(F_i)$ then the identity for $M_1(X)$ is given by
$$M_1(X) = \Vol(T_1(X)) = \sum_i \Omega(E_i).$$
In the paper \cite{BK10}, we show that this is the identity  in equation \ref{bridgeman}.  This proves  part 3) of the Main Theorem.

\section{Zero Moment identity is Basmajian's identity}
As $M_0(X) = \mu(G(X))$ the identity for $M_0(X)$ is
$$\mu(G(X)) = \sum_i \mu(F_i)$$
We now prove some properties of the Liouville measure $\mu$ needed to evaluate both sides of this identity and show that it gives Basmajian's identity in equation \ref{basmajian}.

\begin{lemma}
Let $\mu$ be the Liouville measure on $G(\Hyp^n)$ and $P$ a plane in $\Hyp^n$. For $U \subseteq P$,  let $G(U) \subseteq G(\Hyp^n)$ be the set of geodesics intersecting $U$ transversely. Then the measure $\mu_P$ on $P$ defined by $\mu_P(U)  = \mu(G(U))$ is a constant times area measure on $P$. In particular there is a constant $K_n > 0$ depending only on dimension such that $\mu_P(U) = K_n.\Vol(U)$.
\label{liouville}
\end{lemma}

{\bf Proof:} In the case of $n=2$ this is a standard property of $\mu$ (see \cite{Bon88}). In general we see that $\mu_P$ gives a Mobius  invariant measure on the hyperbolic plane $P$ and therefore must be a multiple of area measure.  Thus $\mu_P (U)  = K_n. \Vol(U)$ where $K_n$ only depends on the dimension $n$.
We calculate $K_n$ in a later section
\eproof

\begin{figure}[htbp] 
   \centering
   \includegraphics[width=5in]{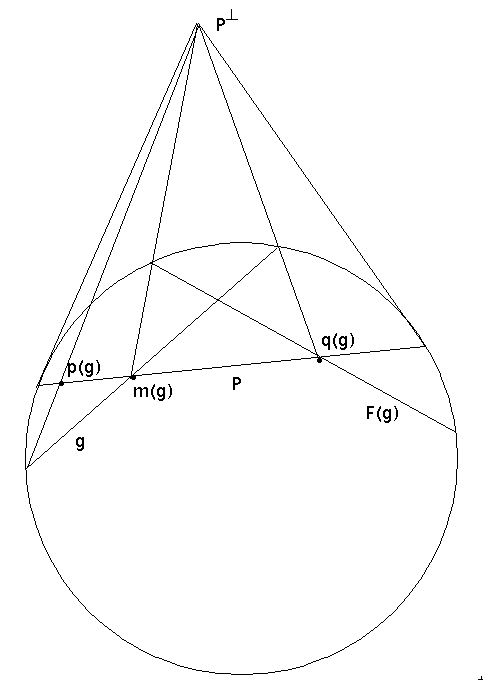}
   \caption{The involution $F$ in the Klein model}
   \label{involution}
\end{figure}

{\bf A Liouville measure preserving map:} We let  $P$ be a plane in $\Hyp^n$ and $G_P$ the set of geodesics intersecting $P$ transversely. Then $G_P =(D_-\times D_+) \cup (D_+\times D_-)$ where $D_-,D_+$ are disjoint open disks in $\Sph_\infty^{n-1}$ with $\partial D_- = \partial D_+$. We let  $g = (x(g),y(g)) \in G_P$ be the ordered pair of endpoints of $g$ and $m(g) = g \cap P$ the point of intersection with $P$. We project the endpoints of $g$ orthogonally onto $P$ to obtain the ordered pair of points $(p(g), q(g)) \in P\times P$. We note that $m(g)$ is the midpoint of the hyperbolic geodesic arc in $P$ joining $p(g),q(g)$. We now define a map $F:G_P \rightarrow G_P$ as follows; We first define $F:D_-\times D_+ \rightarrow D_-\times D_+$ and extend it to $D_+\times D_-$
by conjugating with the map switching endpoints, i.e. if $i(x,y) = (y,x)$ then $F(g) = i(F(i(g))$ for $g \in D_+\times D_-$. If $g \in D_-\times D_+$ then we define $F(g) = h$ where $h$ is the unique geodesic (in $D_-\times D_+$) such that $m(h) = q(g)$ and $q(h) = m(g)$ (see figure \ref{involution}). By construction, $F$ is a involution and if $u$ is an isometry of $\Hyp^n$ fixing $P$, then $F$ commutes with $u$. We note that $F$ is not the action of an isometry on the space of geodesics.

\begin{lemma}
The map $F:G_P \rightarrow G_P$ preserves Liouville measure.
\end{lemma}

{\bf Proof:}
We show that $F_*(\mu) = \mu$ by showing the Radon-Nikodym derivative
$$\nu(x) = \frac{dF_*(\mu)}{d\mu} (x) = 1.$$

Alternately, $\nu(x)$ is the function such that for any $\phi$ smooth compactly supported function on $G_P$ then
$$F_*(\mu)(\phi) = \int \phi(x) d(F_*(\mu)) =  \int \phi(x) \nu(x) d\mu(x)$$
Therefore
$$F_*(\mu)(\phi) = \mu(\phi \circ F) = \int \phi(F(x))d\mu(x) = \int \phi(x) \nu(x) d\mu(x)$$

As $F$ is an involution $F(F(x)) = 1$,  then by the change of variables formula
we have
$$\int \phi(x) d\mu(x) = \int \phi(F(F(x))d\mu(x) = \int \phi(F(x))\nu(x)d\mu(x) = \int \phi(x)\nu(F^{-1}(x))\nu(x)d\mu(x)$$ Therefore
$$ \nu(x)\nu(F^{-1}(x)) = 1 \qquad \mbox{or}\qquad \nu(x)\nu(F(x)) = 1.$$

Similarly we have if  $u$ is a hyperbolic isometry  then $u$ preserves the Liouville measure and $u_*(\mu) = \mu$. If $u$ also fixes $P$ then $F$ commutes with $u$ and $F(u(x)) = u(F(x))$. Therefore by the change of variables formula again we have
$$\nu((u(x)) = \nu(x).$$

Thus combining the above, if there exists an isometry $u$ fixing $P$ such that $u(x) = F(x)$, then $\nu(x) = \nu(F(x))$. But as $\nu(x)\nu(F(x)) = 1$ we then obtain $\nu(x) = 1$.

To find such an isometry, we note that geodesic for $x \in G_P$, and $y = F(x)$ we have the four points $p(x), m(x) = q(y), q(x) = m(y), p(y)$ all collinear on $P$. We choose the plane $P^*$ perpendicular to $P$ which bisects the hyperbolic interval $[p(x),p(y)]$. Then refection in  $P^*$ fixes $P$ and sends $x$ to $y = F(x)$.

Thus $\nu(x) = 1$ for all $x \in G_P$ and therefore $F_*(\mu) = \mu$.
\eproof

\begin{corollary}
If $\alpha$ is an orthogeodesic  then
$$\mu(F_\alpha) = K_n.V_{n-1}\left(\log\left(\coth{\frac{l(\alpha)}{2}}\right)\right)$$
\label{equiv}
\end{corollary}

{\bf Proof:}
We consider disjoint planes $P,Q$ with perpendicular distance equal the ortholength $l(\alpha)$. Then $\mu(F_\alpha) = \mu(S)$ where $S \subseteq G_P$ of geodesics which intersect $Q$. We let $G_P =(D_-\times D_+) \cup (D_+\times D_-)$ where $\partial Q \subseteq D_+$.
 By the above lemma, $\mu(S) = \mu(F(S))$. Let $B$ be the orthogonal projection of $Q$ onto $P$. Then by elementary hyperbolic geometry, $B$ is a ball of radius $r = \log(\coth{\frac{l(\alpha)}{2}})$.
The set  $F(S)$ is precisely the set of geodesics in $G_P$ transversely intersecting $B$ which we denote by $G(B)$. To see this, we note that if $g \in S$  then $F(g)$ intersects $P$ in $B$ giving $F(S) \subseteq G(B)$. Similarly if $g \in G(B)$ then $F(g)$ is in $S$. Thus, as $F$ is an involution $F(S) = G(B)$.

Therefore by lemma \ref{liouville}
$$\mu(F_\alpha) = \mu(S) = \mu(F(S)) = \mu(G(B)) = K_n\Vol(B) = K_nV_{n-1}\left(\log\left(\coth{\frac{l(\alpha)}{2}}\right)\right)$$
\eproof

We will now prove the second part of the Main theorem.

By the above corollary \ref{equiv} we obtain
$$\mu(G(X)) =  \sum_{\alpha \in O_X}\mu(F_\alpha) =  K_n\sum_{\alpha \in O_X}V_{n-1}\left(\log\left(\coth{\frac{l(\alpha)}{2}}\right)\right)$$

For each boundary component $B_i$ of $\partial X$ we let $B'_i$ in $\Hyp^n$ be a hyperplane which is a lift of $B_i$. We further take a fundamental domain $D_i$ on $B'_i$ for the action of $\Gamma$. We let $C_i$ be the set of geodesics which intersect $D_i$ transversely such that the geodesics point into $\tilde{X}$ on $B'_i$. Let $G = \cup C_i$ then we see that $G$ is a lift of $G(X)$ (except for a set of measure zero). To see this note that for almost every $g \in G(X)$, $g$ is a proper arc from a $B_i$ to a $B_j$ where the orientation of $g$ is pointing into $X$ at $B_i$ and out at $B_j$. Therefore $g$ has lift $g'$ in $C_i$.  Also for $i\neq j$ if  $g \in C_i \cap C_j$ then $g$ points inward on both $B_i$ and $B_j$. Therefore $C_i \cap C_j = \emptyset$ for $i \neq j$. Also if $g_1, g_2 \in C_i$ are lifts of the same element of $G(X)$ then there is a $\gamma \in \Gamma$ with $g_2 = \gamma(g_1)$. As $D_i$ is a fundamental domain for the action of $\Gamma$ on $B'_i$ then $g_1, g_2$ must have endpoints on the boundary of $D_i$ which is measure zero.

Using this we can calculate $\mu(G(X))$. We have
$$\mu(G(X)) = \sum_i\mu(C_i)$$
By the above lemma \ref{liouville}, $\mu(C_i) = (K_n/2)\Vol(B_i)$  where  the factor of two comes from $C_i$ containing half the geodesics in the set $G(B_i)$ (those pointing into $\tilde{X}$). Therefore
$$ \mu(G(X)) = \frac{K_n}{2} \sum \Vol(B_i) = \frac{K_n}{2}\Vol(\partial X)$$
Thus
$$ \frac{K_n}{2}\Vol(\partial X) = K_n\sum_{\alpha \in O_X}V_{n-1}\left(\log\left(\coth{\frac{l(\alpha)}{2}}\right)\right)$$
giving Basmajian's identity.
\eproof

  \section{Calculating $K_n$}
 To calculate the constant $K_n$, we derive the Lebesgue density of $\mu_P$ at a point $p$. Let $p$ be at the origin of the Poincare model of $\Hyp^n$ and let $P$ be the horizontal plane through $p$. Let $B = B_{n-1}(p,r)$ be a small $(n-1)$-dimensional ball in $P$ about $p$. Then we have
$$d\mu = \frac{2dV_xdV_y}{|x-y|^{2n-2}}$$
where $V$ is volume measure on the unit sphere. If we consider geodesics $G(B)$ then if $(x,y) \in G(B)$, we have $|x-y| \simeq 2$. Also if we let $\theta_x$ be the angle the ray $px$ from $x$ to $p$ makes with the plane $P$, then the set $G(B)_x = \{y \in \Sph^{n-1} \ | \ (x,y) \in G(B)\}$ is a small ball in the unit sphere. In fact the set $G(B)$ is the image of $B$ under stereographic projection from $x$. Therefore $G(B)_x$ is an $(n-1)$-dimensional ellipsoid with the axes of $B$ perpendicular to ray $px$ being approximately $2r$ and other axis  approximately $2r\sin(\theta_x)$. Therefore
$$\Vol(G(B_x)) \simeq 2^{n-1}\sin(\theta_x).\Vol(B^s_{n-1}(r))$$ where $B^s_{k}(r)$ is a k-dimensional ball of radius $r$ in the unit sphere.
Integrating we get
$$\mu(G(B)) = \int_{\Sph^{n-1}} 2.\left( \int_{G(B)_x} \frac{dV_y}{|x-y|^{2n-2}} \right)dV_x\simeq  \frac{Vol(B^s_{n-1}(r))}{2^{n-2}} \int_{\Sph^{n-1}} \sin(\theta_x) dV_x$$

Thus if $A$ is area measure on $P$ then
$$K_n = \frac{d\mu_P}{dA} = \frac{1}{2^{n-2}} \int_{\Sph^{n-1}} \sin(\theta_x)dV_x$$
The set $S_t = \{x\ | \ \theta_x = t\}$ is an $(n-2)$-dimensional sphere of radius $|\cos(t)|$. Therefore  as $dV_x = d\theta dV_{S_t}$
 $$K_n = \frac{\Vol(\Sph^{n-2})}{2^{n-2}} \int_{0}^{\pi} \sin(\theta)|\cos^{n-2}(\theta)| d\theta = \frac{\Vol(\Sph^{n-2})}{2^{n-1}} \int_{0}^{\pi/2} \sin(\theta)\cos^{n-2}(\theta)d\theta =  \frac{\Vol(\Sph^{n-2})}{2^{n-1}(n-1)} $$
In terms of the Gamma function  $\Gamma$ we have
$$ \Vol(\Sph^k) =  \frac{2\pi^{\frac{k+1}{2}}}{\Gamma(\frac{k+1}{2})} $$ Giving
$$K_n = \frac{\pi^{\frac{n-1}{2}}}{2^{n-1}\Gamma(\frac{n+1}{2})}$$
Note: We have $K_2 = 1$.

\section{Explicit Integral Formulae for $F_{n,k}$} In \cite{B11} we derive a formula for $L$ in the surface case. Using this we can write
$$F_{2,k}(l) = \frac{1}{2^{k-2}} \int_0^{a}\int_1^\infty \frac{\log\left|\frac{y(y-a)(x-1)}{x(x-a)(y-1)}\right|^{k}}{(y-x)^2}dxdy.$$
where $a = \sech^2(l/2)$.

In \cite{BK10} we consider the case $n > 2$ where we derive the explicit formula for $L$ and reduce the integral formula to a triple integral via an elementary substitution.  Using this we can also reduce the integral of $F_{n,k}$ to a triple integral of the form
$$F_{n,k}(l) = \frac{\Vol(\Sph^{n-2})\Vol(\Sph^{n-3})}{2^{k-2}}\int_0^1 \frac{r^{n-3}}{(1-r^2)^{\frac{n-2}{2}}}dr\int^{1}_{-1} du\int^{\infty}_{b}  \frac{\log\left(\frac{(v^2-1)(u^2-b^2)}{(v^2-b^2)(u^2-1)}\right)^{k}}{(v-u)^{n}}dv$$
where $b = \sqrt{\frac{e^{2l}-r^2}{1-r^2}}$.

\section{Moment Generating function}
The moment generating function of a random variable $Y$ is the function $f_Y(t) =E[e^{tY}]$ where $E$ is the expected value. We define the moment generating function for measure $M$
$$\zeta^M_X(t) = M(e^{tx}) =  \int_0^\infty e^{tx} dM = \int_{G(X)} e^{t.L(g)} \ d\mu.$$
It follows from above that
$$\zeta^M_X(t) = \sum_{l \in L_X} F_n(t,l)$$
for some function $F_n$ depending only on the dimension $n$.
We similarly can define $\zeta^N_X(t) = N(e^{xt})$. The it follows that the two functions are related by
$$\zeta_X^N(t) = \frac{d}{dt}\left(\zeta_X^M(t)\right) .$$

{\bf A simple example, the ideal triangle:} We consider the case of $X$ being an ideal triangle. In this case, it is more natural to consider measure $N$ (in particular $M$ has infinite mass). In \cite{BD07} we show that
$$dN= \frac{12x^2}{\sinh^2{x}}dx.$$
and that
$$N_k(X) =  \frac{3(k+2)\\!\zeta(k+2)}{2^{k-1}}$$
where $\zeta$ is the Riemann zeta function.  In particular the average time to the boundary is
 $$A(X) = \frac{N_1(X)}{2\Vol(T_1(X))} = \frac{9}{2\pi^2}\zeta(3).$$

It follows by integrating that
$$\zeta^N_X(t) =  \int_0^\infty\frac{12x^2e^{xt}}{\sinh^2{x}}dx = 12 \left(\zeta(2,1 - \frac{t}{2}) + \frac{t}{2} \zeta(3,1 - \frac{t}{2})\right)$$
where  $\zeta(s,t)$ is the Hurwitz  zeta function
$$\zeta(s,t) = \sum_{k=0}^\infty \frac{1}{(k+t)^s}$$

\section{The Surface Case}

For the surface case, the identities in the Main Theorem can be written in terms of polylogarithm functions.

{\bf Polylogarithms:} The $k^{th}$ polylogarithm function $\mbox{Li}_k$ is defined by the Taylor series
$$\mbox{Li}_{k}(z) = \sum_{i=1}^{\infty} \frac{z^n}{n^k }$$
for $|z| < 1$ and by analytic continuation to $\Cplex$. In particular
$$\mbox{Li}_{0}(z) = \frac{z}{1-z} \qquad \qquad
\mbox{Li}_{1}(z) = -\log(1-z).$$
Also
$$\mbox{Li}'_{k}(z) = \frac{\mbox{Li}_{k-1}(z)}{z} \qquad \mbox{giving} \qquad \mbox{Li}_{k}(z) = \int_{0}^{z} \frac{\mbox{Li}_{k-1}(t)}{t}\ t.$$
Also the functions $\mbox{Li}_{k}$ are related to the Riemann $\zeta$ function by $\mbox{Li}_{k}(1) = \zeta(k)$.

Below  we describe some properties of the dilogarithm and trilogarithm function. They can all be found in the 1991 survey "Structural Properties of Polylogarithms" by L. Lewin (see \cite{Lew91}).

{\bf Dilogarithm:}
The dilogarithm function $\mbox{Li}_2(z)$ is given by
$$\mbox{Li}_{2}(z) = -\int_{0}^{z} \frac{\log(1-t)}{t}\ dt.$$
 From the power series representation, it is easy to see that the dilogarithm function satisfies the functional equation
$$\mbox{Li}_{2}(z) + \mbox{Li}_{2}(-z) = \frac{1}{2}\mbox{Li}_{2}(z^{2}).$$
Other functional relations of the dilogarithm can be best described by normalizing the dilogarithm function.
The (extended) Rogers dilogarithm function (see \cite{Rog07}) is defined by
$$R(x) = \mbox{Li}_{2}(x) + \frac{1}{2}\log|x|\log (1-x)\qquad x \leq 1.$$
This function arises in calculating hyperbolic volume as the imaginary part of $R(z)$ is the volume of the hyperbolic tetrahedron with vertices having cross ratio $z$.

Also in terms of the Rogers function, various identities have nice form.
Euler's reflection relations  for the dilogarithm are given by
\begin{eqnarray}
R(x) + R(1-x) = R(1) = \frac{\pi^2}{6} \qquad 0 \leq x \leq 1 \nonumber \\
 R(-x) + R(-x^{-1}) = 2R(-1) = -\frac{\pi^{2}}{6}  \qquad x >  0
 \label{euler}
\end{eqnarray}
Also Landen's identity is
\begin{equation}
R\left(\frac{-x}{1-x}\right) = -R(x)  \qquad  0 < x < 1
  \label{landen}
  \end{equation}
  and Abel's functional equation  is
  \begin{equation}
  R(x) + R(y) = R(xy) + R\left(\frac{x(1-y)}{1-xy}\right) + R\left(\frac{y(1-x)}{1-xy}\right).
  \label{abel}
  \end{equation}

In \cite{B11}, we showed that the orthospectra of a hyperbolic surface satisfies the following generalized orthospectrum identity.
\begin{theorem}{(Bridgeman, \cite{B11})}
Let $S$ be a finite area hyperbolic surface with non-empty totally geodesic boundary and $C_S$ boundary cusps. Then
$$ \sum_{l \in L_S} R\left(\frac{1}{\cosh^2\frac{l}{2}}\right) = \frac{\pi^{2}(6|\chi(S)| -C_{S})}{12}$$
\label{Bidentity}
\end{theorem}

{\bf Trilogarithm:}

By definition, the trilogarithm function is given by
$$Li_3(z) = \int_0^z \frac{Li_2(t)}{t}\ dt$$

The trilogarithm also satisfies a number of identities.
\begin{equation}
Li_3(z) + Li_3(-z) =  \frac{1}{4}Li_3(z^2)
\label{tri1}
\end{equation}
\begin{equation}
Li_3(-z) - Li_3(-z^{-1}) = - \frac{1}{6}\log^3(z) - \frac{\pi^2}{6}\log(z)
\label{tri2}
\end{equation}
and
\begin{equation}
Li_3(z)+Li_3(1-z) - Li_3(1-z^{-1}) = \zeta(3)+ \frac{1}{6}\log^3(z) + \frac{\pi^2}{6}\log(z)-\frac{1}{2}\log^2(z)\log(1-z)
\label{tri3}
\end{equation}

If $l_i$ is an ortholength of $S$, we define $$a_i =\frac{1}{\cosh^2\left(\frac{l_i}{2}\right)}.$$
We will often use the spectrum $\{a_i\}$ instead of $\{l_i\}$. In the paper \cite{B11}, we studied the measure $N$ and derived the following;
\begin{theorem}{(Bridgeman, \cite{B11})}
There exists a smooth function $\rho:\Real_+\times(0,1) \rightarrow\Real_+$ such that
$$dN = \rho_S(x)dx = \left(\frac{4C_{S}x^{2}}{\sinh^{2}(x)} + \sum_{a_i}\rho(x,a_i)\right)  dx.$$
Furthermore
$$\int \phi(x) \rho(x,a)dx = \int_0^{a}\int_1^\infty \frac{4\phi(L_a(x,y))L_a(x,y)}{(y-x)^2}dxdy$$
where
$$L_a(x,y) = \frac{1}{2}\log\left|\frac{y(y-a)(x-1)}{x(x-a)(y-1)}\right|$$
\end{theorem}

{\bf The Moment Identities for Surface:}
If we apply the above theorem to the function $\phi(x) = 1$ then we recover the identity in theorem \ref{Bidentity}.

To find $M_k(S)$,  we note that $M_{k}(S) = N_{k-1}(S)$.  Therefore by the above, we let $\phi(x) = x^{k-1}$ to get
$$M_k(S) = N_{k-1}(S) =  \int_0^\infty x^{k-1}\rho_S(x)dx = \sum_{a_i} \left(\int_0^{a_i}\int_1^\infty \frac{4L^{k}_{a_i}(x,y)}{(y-x)^2}dxdy\right)+4C_S\left(\int_0^\infty \frac{x^{k+1}dx}{\sinh^2(x)}\right)$$
We then define
$$F_k(a) = \int_0^{a}\int_1^\infty \frac{4L_{a}(x,y)^{k}}{(y-x)^2}dxdy = \frac{1}{2^{k-2}}\int_0^{a}\int_1^\infty \frac{\log\left|\frac{y(y-a)(x-1)}{x(x-a)(y-1)}\right|^{k}}{(y-x)^{2}}dxdy.$$
Integrating  we have
$$\int_0^\infty \frac{x^{k+1}dx}{\sinh^2(x)} = \frac{(k+1)\\!\zeta(k+1)}{2^{k}}$$
Therefore we obtain
$$M_k(S) = \left(\sum_{l\in L_S} F_k\left(\frac{1}{\cosh^2(\frac{l}{2})}\right) + C_S \frac{(k+1)\\!\zeta(k+1)}{2^{k-2}}\right).$$
We note that as $\zeta(1) = \infty$, if $S$ has boundary cusps then $M_0(S)$ is not finite. This corresponds to the fact that $M_0(S) = Length(\partial S)$ which is infinite in the case of boundary cusps.

Functions $F_0, F_1$ are given in terms of simple logarithms and dilogarithms respectively. An induction argument shows that $F_k$ can be written as a sum of polylogarithm functions of order at most $k+1$. We will calculate an explicit formula for $F_2$ in terms of trilogarithms  in the next section. This will give us the formula for the average hitting time $A(S)$ for geodesic flow described in theorem \ref{avehit}.

\section{A Somewhat Brutal Calculation}
We will now obtain an explicit formula for the average hitting time in the surface case in terms of sums of polylogarithms evaluated at ortholengths. We let $F = F_2$ given by the above integral formula.  Then
$$F(a) =  \int_0^{a}\int_1^\infty \frac{\log\left|\frac{y(y-a)(x-1)}{x(x-a)(y-1)}\right|^{2}}{(y-x)^{2}}dxdy.$$
 Using Mathematica to calculate the indefinite integral first, gives 7,858 polylogarithm terms  which then need to be evaluated at the 4 limits to give a final total of approximately 30,000 terms. Also the terms must be grouped so that evaluation gives a finite limit.  As this seems a daunting task, we do the calculation directly using hyperbolic relations to simplify as we go along. The calculation is somewhat tedious but the final answer surprisingly short.

\begin{figure}[htbp] 
   \centering
   \includegraphics[width=5in]{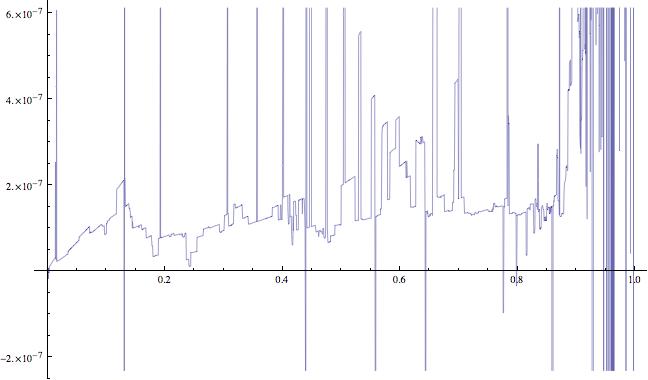}
   \caption{Difference Between Numerical Integration of F and Polylogarithm Formula}
   \label{numerical}
\end{figure}

For the reader who would rather skip the long and tedious calculation, evidence for its validity is given by   figure \ref{numerical}, which  is a plot of the difference between the polylogarithm formula for $F(a)$ and  and its values using numerical integration.  As can be seen from the plot, the difference is less than $10^{-6}$ indicating they are the same function.

We have that for $a \in (0,1)$
$$F(a)  = \int_0^a \int_1^\infty \frac{\log\left(\frac{y(y-a)(y-1)}{x(x-a)(y-1)}\right)^2}{(y-x)^2}dydx$$
Decomposing into cross-ratios, we have
$$ F(a) = \int_0^a \int_1^\infty \frac{\left(\log\left|\frac{y(x-1)}{x(y-1)}\right|+\log\left|\frac{y-a}{x-a}\right|\right)^2}{(y-x)^2}dydx =\int_0^a \int_1^\infty \frac{\log\left|\frac{y(x-1)}{x(y-1)}\right|^2 +2.\log\left|\frac{y(x-1)}{x(y-1)}\right|.\log\left|\frac{y-a}{x-a}\right| +\log\left|\frac{y-a}{x-a}\right|^2}{(y-x)^2}dydx$$
Under the mobius transformation $m(z) = a/z$, we let $X= m(x), Y = m(y)$, then by invariance of cross ratios,
$$\frac{y(x-1)}{x(y-1)} = \frac{(y-0)(x-1)}{(x-0)(y-1)} = \frac{(m(y)-m(0))(m(x)-m(1)}{(m(x)-m(0)(m(y)-m(1)} = \frac{(Y-\infty)(X-a)}{(X-\infty)(Y-a)} = \frac{X-a}{Y-a}$$
Thus
$$\int_0^a \int_1^\infty \frac{\log\left|\frac{y(x-1)}{x(y-1)}\right|^2}{(y-x)^2}dydx = \int_\infty^1 \int_a^0 \frac{\log\left|\frac{X-a}{Y-a}\right|^2}{(Y-X)^2}dYdX = \int_0^a \int_1^\infty \frac{\log\left|\frac{y-a}{x-a}\right|^2}{(y-x)^2}dydx$$
$$ F(a) =   2.\int_0^a \int_1^\infty \frac{\log\left|\frac{y-a}{x-a}\right|^2}{(y-x)^2}dydx + 2.\int_0^a \int_1^\infty \frac{\log\left|\frac{y(x-1)}{x(y-1)}\right|.\log\left|\frac{y-a}{x-a}\right|}{(y-x)^2}dydx $$
We write this as $F(a) = 2I_1(a) + 2I_2(a)$, where $I_1, I_2$ are the above integrals.

In order to calculate the above integrals we will need the following integral equations.
\begin{eqnarray}
\int \frac{\log(x)}{x-a} dx  &=& \log(x)\log(1-x/a)+Li_2(x/a)\\
\int \frac{\log(x)^2}{x-a} dx &=& \log(x)^2\log(1-x/a)+2\log(x)Li_2(x/a)-2Li_3(x/a)\\
\int \frac{\log(x)\log(x-a)}{x} dx &=& \frac{1}{2}\log(x)^2\log(a)-\log(x)Li_2(x/a)+Li_3(x/a)
\end{eqnarray}

\subsection{Integral $I_1$}
\begin{lemma}
$$I_1(a) = -\log(1-a)\log^2(a)+\log^3(1-a)-4\log\left(\frac{a}{1-a}\right)Li_2(a)-6Li_3\left(\frac{a}{a-1}\right)$$
\end{lemma}
{\bf Proof:}
We decompose $I_1$ to obtain
$$I_1 = \int_0^a\int_1^\infty \frac{ \log(y-a)^2 - 2\log(y-a)\log(a-x) + \log(a-x)^2}{(y-x)^2}dydx = J_1 -2J_2 + J_3$$
{\bf Integral $J_1$:}\newline
We have
$$J_1 = \int_1^\infty  \log(y-a)^2\left(\int_0^a \frac{dx}{(y-x)^2}\right)dy  = \int_1^\infty \log(y-a)^2\left(\frac{1}{y-a} - \frac{1}{y}\right)dy.$$
Then
$$ \int \frac{\log(y-a)^2}{y-a}dy = \frac{\log(y-a)^3}{3}$$
Also by integral equations above
$$\int \frac{\log(y-a)^2}{y}dy = \log\left(\frac{y}{a}\right)\log(y-a)^2+2\log(y-a)Li_2\left(1-\frac{y}{a}\right)-2Li_3\left(1-\frac{y}{a}\right)$$
Thus
$$J_1 = -\frac{\log(1-a)^3}{3}+  \frac{\pi^2}{3}\log(a)+\frac{1}{3}\log(a)^3+
\log\left(\frac{1}{a}\right)\log(1-a)^2+2\log(1-a)Li_2\left(1-\frac{1}{a}\right)-2Li_3\left(1-\frac{1}{a}\right)$$

{\bf Integral $J_2$:}
By parts we have
$$J_2 = \int_0^a \log(a-x)\left(\int_1^\infty \frac{\log(y-a)}{(y-x)^2}dy\right)dx = \int_0^a\log(a-x)\left(\frac{\log(1-a)}{1-x}-\frac{\log(1-a)}{a-x}+\frac{\log(1-x)}{a-x}\right)dx$$
As above  we have
$$\int _0^a \frac{\log(a-x)}{1-x} dx = \left.\left(-\log(a-x)\log\left(\frac{1-x}{1-a}\right) -Li_2\left(\frac{a-x}{a-1}\right)\right)\right|_0^a = -\log(a)\log(1-a)+Li_2\left(\frac{a}{a-1}\right).$$
Also
$$\int\frac{\log(a-x)}{a-x}dx = -\frac{1}{2}\log(a-x)^2$$
$$\int \frac{\log(a-x)\log(1-x)}{a-x}dx =-\frac{1}{2}\log(a-x)^2\log(1-a)+\log(a-x)Li_2\left(\frac{a-x}{a-1}\right) - Li_3\left(\frac{a-x}{a-1}\right) $$
Combining we get
$$\int_0^a\frac{\log(a-x)(\log(1-x)-\log(1-a))}{a-x}dx =\left. \log(a-x)Li_2\left(\frac{a-x}{a-1}\right) - Li_3\left(\frac{a-x}{a-1}\right) \right|_0^a = -\log(a)Li_2\left(\frac{a}{a-1}\right)+Li_3\left(\frac{a}{a-1}\right) $$
Thus
$$J_2 = -\log(a)\log(1-a)^2-\log\left(\frac{a}{1-a}\right)Li_2\left(\frac{a}{a-1}\right)+Li_3\left(\frac{a}{a-1}\right)$$
{\bf Integral $J_3$:}

Similarly to $J_1$ we have
$$J_3 = \int_0^a \log(a-x)^2 \left(\int_1^\infty\frac{1}{(y-x)^2}dy\right) dx =  \int_0^a \frac{\log(a-x)^2}{1-x} dx$$
Therefore
\begin{eqnarray*}
J_3 =\left. -\log(a-x)^2\log\left(\frac{1-x}{1-a}\right) - 2\log(a-x)Li_2\left(\frac{a-x}{a-1}\right)+2Li_3\left(\frac{a-x}{a-1}\right)\right|_0^a\\
-\log(1-a)\log(a)^2+2\log(a)Li_2\left(\frac{a}{a-1}\right)-2Li_3\left(\frac{a}{a-1}\right)
\end{eqnarray*}

Thus we have the formula for $I_1 = J_1-2J_2+J_3$ giving
\begin{eqnarray*}
I_1 &=  \left(-\frac{\log(1-a)^3}{3}+ \frac{\pi^2}{3}\log(a) +\frac{1}{3}\log(a)^3-\log(a)\log(1-a)^2+2\log(1-a)Li_2\left(1-\frac{1}{a}\right)-2Li_3\left(1-\frac{1}{a}\right)\right)\\ &-2\left(-\log(a)\log(1-a)^2-\log\left(\frac{a}{1-a}\right)Li_2\left(\frac{a}{a-1}\right)+Li_3\left(\frac{a}{a-1}\right)\right)\\&+\left(-\log(1-a)\log(a)^2+2\log(a)Li_2\left(\frac{a}{a-1}\right)-2Li_3\left(\frac{a}{a-1}\right)\right)
\end{eqnarray*}
Using the identities \ref{euler}, \ref{tri2}, we get
$$I_1(a) =  -\log(1-a)\log\left(\frac{a}{1-a}\right)^2+4\log\left(\frac{a}{1-a}\right)Li_2\left(\frac{a}{a-1}\right)-6Li_3\left(\frac{a}{a-1}\right)$$
Simplifying further we also get
$$I_1(a) = -\log(1-a)\log^2(a)+\log^3(1-a)-4\log\left(\frac{a}{1-a}\right)Li_2(a)-6Li_3\left(\frac{a}{a-1}\right)$$
\eproof
\subsection{Integral $I_2$}
\begin{lemma}
\begin{eqnarray*}
I_2(a) &=&2\zeta(3)+\frac{2\pi^2}{3}\log(1-a)+\frac{1}{3}\log^3(1-a)- \log^2(1-a)\log(a)-\log^2(a)\log(1-a)\\
&&-4\log(a)Li_2(a)+4Li_3(a)-2Li_3\left(\frac{-a}{1-a}\right)-2Li_3(1-a)
\end{eqnarray*}
\end{lemma}

{\bf Proof}
We have
$$I_2 = \int_0^a \int_1^\infty\frac{\log\left|\frac{y(x-1)}{x(y-1)}\right|.\log\left|\frac{y-a}{x-a}\right|}{(y-x)^2}dydx$$

Again we decompose into integrals
$$J_1 =  \int_0^a \int_1^\infty \frac{\log\left(\frac{y}{y-1}\right)\log(y-a)}{(y-x)^2}dydx$$
$$J_2 =  \int_0^a \int_1^\infty \frac{\log\left(\frac{1-x}{x}\right)\log(y-a)}{(y-x)^2}dydx$$
$$J_3 =  \int_0^a \int_1^\infty \frac{\log\left(\frac{x}{1-x}\right)\log(a-x)}{(y-x)^2}dydx$$
$$J_4 =  \int_0^a \int_1^\infty \frac{\log\left(\frac{y-1}{y}\right)\log(a-x)}{(y-x)^2}dydx$$

As before we can integrate to get
$$J_1 = \int_1^\infty \log\left(\frac{y}{y-1}\right)\log(y-a)\left(\frac{1}{y-a}-\frac{1}{y}\right)dy$$
and
$$J_3 =  \int_0^a \frac{\log\left(\frac{x}{1-x}\right)\log(a-x)}{1-x}dx$$
We integrate in $J_2, J_3$ the simple factor to get

$$J_2 =  \int_0^a \log\left(\frac{1-x}{x}\right)\left(\int_1^\infty \frac{\log(y-a)}{(y-x)^2}dy\right)dx = \int_0^a \log\left(\frac{1-x}{x}\right)\left(\frac{\log(1-x)-\log(1-a)}{a-x}+\frac{\log(1-a)}{1-x}\right)dx$$
Combining $J_2, J_3$ we get
$$J_2 + J_3 = \int_0^a \log\left(\frac{1-x}{x}\right)\left(\frac{\log\left(\frac{1-x}{1-a}\right)}{a-x}-\frac{\log\left(\frac{a-x}{1-a}\right)}{1-x}\right)dx$$

The Roger's normalized dilogarithm is given by $R(x) = Li_2(x) + \frac{1}{2}\log|x|\log(1-x)$ for $x <1$. We note that
$$R'(x) = -\frac{1}{2}\left(\frac{\log(1-x)}{x}+\frac{\log|x|}{1-x}\right)$$
Thus we let $R_1(x) = R(\frac{a-x}{1-x})$ and note that
$$R_1'(x) = \frac{1}{2}\left(\frac{\log\left(\frac{a-x}{1-a}\right)}{1-x}-\frac{\log\left(\frac{1-x}{1-a}\right)}{a-x}\right)$$
Thus
$$J_2 + J_3 = 2\int_0^a \log\left(\frac{1-x}{x}\right)R_1'(x) dx$$
Considering $J_4$ we have similarly
$$J_4 = \int_1^\infty   \log\left(\frac{y-1}{y}\right) \left(\int_0^a \frac{(\log(a-x)}{(y-x)^2}dx\right)dy =\int_1^\infty  \log\left(\frac{y-1}{y}\right) \left(\frac{\log(y-a)-\log(y)+\log(a)}{y-a}-\frac{\log(a)}{y}\right)dy  $$
We then have
$$J_1+J_4 = \int_1^\infty  \log\left(\frac{y-1}{y}\right) \left(\frac{\log\left(\frac{a}{y}\right)}{y-a}+\frac{\log\left(\frac{y-a}{a}\right)}{y}\right)dy  $$
We now let $R_2(y) = R(a/y)$, then
$$R_2'(y) =   \frac{1}{2}\left(\frac{\log\left(\frac{a}{y}\right)}{y-a}+\frac{\log\left(\frac{y-a}{a}\right)}{y}\right)$$
giving
$$J_1+J_4 = 2\int_1^\infty   \log\left(\frac{y-1}{y}\right)R_2'(y)dy$$
Letting $x = a/y$ we have

$$J_1 + J_4 = 2\int_a^0 \log\left(\frac{a-x}{a}\right)R'(x)dx = 2\log(a)R(a)- 2\int_0^a\log(a-x)R'(x)dx  $$

Similarly let $u = (a-x)/(1-x)$ then $x = (a-u)/(1-u)$ and
$$J_2 + J_3 = 2\int_{0}^a \log\left(\frac{1-a}{a-u}\right)R'(u) du = 2\log(1-a)R(a)- 2\int_0^a \log(a-x)R'(x)dx$$
Giving
$$I_2 = 2\log(a(1-a))R(a)-4\int_0^a \log(a-x)R'(x)dx = 2\log(a(1-a)R(a) + 2\int_0^a \log(a-x) \left(\frac{\log(1-x)}{x}+\frac{\log(x)}{1-x}\right)dx$$

We note the formula
\begin{eqnarray*}
G(x,a)&=&\int \frac{\log(a-x)\log(1-x)}{x}dx =  \\
&&\log(1-x)\log\left(\frac{a-x}{a}\right)\log\left(\frac{x}{a}\right) + \log(a)\log(1-x)\log(x)+ \frac{1}{2}\log(a)\log^2(1-x)\\
&+&\log(a(1-x))Li_2(1-x) +\log\left(\frac{a-x}{a}\right)Li_2\left(\frac{a-x}{a}\right) + \log\left(\frac{a-x}{a(1-x)}\right)\left(Li_2\left(\frac{a-x}{1-x}\right)- Li_2\left(\frac{a-x}{a(1-x)}\right)\right)\\
&-&Li_3(1-x)-Li_3\left(\frac{a-x}{1-x}\right) - Li_3\left(\frac{a-x}{a}\right) + Li_3\left(\frac{a-x}{a(1-x)}\right)
\end{eqnarray*}

Thus taking limits we have
$$
 \int_0^a \frac{\log(a-x)\log(1-x)}{x}dx = $$
$$\zeta(3)-\frac{\pi^2}{6}\log(a)+\frac{1}{2}\log(a)\log^2(1-a)+\log(1-a)\log^2(a)+\log(a(1-a))Li_2(1-a) +Li_3(a)-Li_3(1-a)$$
From the above
\begin{eqnarray*}
H(x,a)&=&\int \frac{\log(a-x)\log(x)}{1-x}dx =  -G(1-a,1-x)= \\
&-&\log(x)\log\left(\frac{a-x}{1-a}\right)\log\left(\frac{1-x}{1-a}\right) - \log(1-a)\log(1-x)\log(x)- \frac{1}{2}\log(1-a)\log^2(x)\\
&-&\log((1-a)x)Li_2(x) -\log\left(\frac{a-x}{1-a}\right)Li_2\left(\frac{x-a}{1-a}\right) \\
&-& \log\left(\frac{a-x}{(1-a)x}\right)\left(Li_2\left(\frac{x-a}{x}\right)- Li_2\left(\frac{x-a}{(1-a)x}\right)\right)\\
&+&Li_3(x)+Li_3\left(\frac{x-a}{x}\right) + Li_3\left(\frac{x-a}{1-a}\right) - Li_3\left(\frac{x-a}{(1-a)x}\right)
\end{eqnarray*}
Thus taking limits we have
$$
 \int_0^a \frac{\log(a-x)\log(x)}{1-x}dx = $$
 $$\frac{\pi^2}{6}\log(1-a)-\frac{1}{3}\log^3(1-a) -\log(1-a)\log^2(a)- \log((1-a)a)Li_2(a) +\log\left(\frac{a}{1-a}\right)Li_2\left(\frac{-a}{1-a}\right)+Li_3(a)-Li_3\left(\frac{-a}{1-a}\right)$$
We now combine to obtain
\begin{eqnarray*}
I_2(a) &=&2\zeta(3)+\frac{2\pi^2}{3}\log(1-a)+\frac{1}{3}\log^3(1-a)- \log^2(1-a)\log(a)-\log^2(a)\log(1-a)\\
&&-4\log(a)Li_2(a)+4Li_3(a)-2Li_3\left(\frac{-a}{1-a}\right)-2Li_3(1-a)
\end{eqnarray*}
\eproof

Finally we combine $I_1, I_2$ to get
\begin{eqnarray*}
 F(a)& = 4\zeta(3) +\frac{4\pi^2}{3}\log(1-a)+\frac{8}{3}\log^3(1-a)-4\log(1-a)\log^2(a)-2\log^2(1-a)\log(a)\\
&-8\log\left(\frac{a^2}{1-a}\right)Li_2(a)+8Li_3(a)-4Li_3(1-a)-16Li_3\left(\frac{-a}{1-a}\right)
\end{eqnarray*}
Using the identity \ref{tri3} we get
\begin{eqnarray*}
F(a) = &&-12\zeta(3)-\frac{4\pi^2}{3}\log(1-a)+6\log^2(1-a)\log(a)-4\log(1-a)\log^2(a)\\
&&-8\log\left(\frac{a^2}{1-a}\right)Li_2(a)+24Li_3(a)+12Li_3(1-a)
\end{eqnarray*}

The function has boundary values $F(0) = 0$ and $F(1) = 12\zeta(3)$ and is maximized at $a = .754493$ with value $17.9804$.

Below is a graph of $F$.
 \begin{figure}[htbp] 
    \centering
    \includegraphics[width=5in]{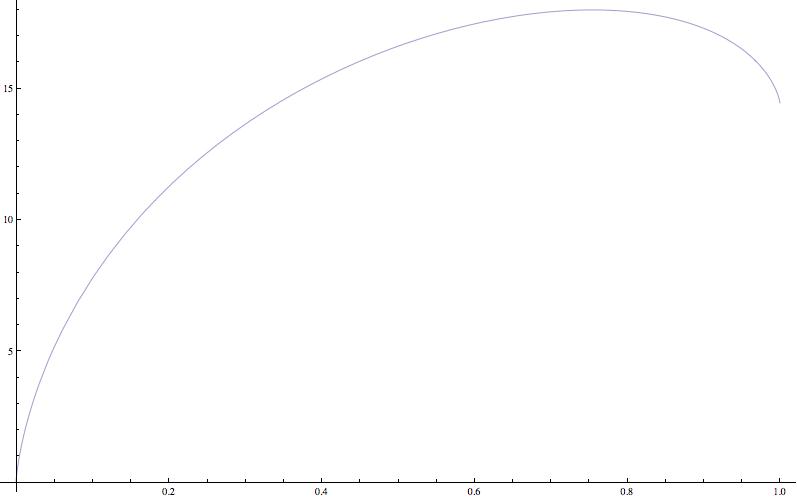}
    \caption{Function F(x)}
    \label{graphF}
 \end{figure}

\end{document}